\documentclass[12pt,reqno]{amsart}
\usepackage{times}
\usepackage{graphicx}
\usepackage{amsmath,amssymb}

\newtheorem{thm}{Theorem}

\newtheorem{lem}[thm]{Lemma}
\newtheorem{prop}[thm]{Proposition}
\newtheorem*{thma}{Theorem A}
\newtheorem*{thmb}{Theorem B}

\newcommand{\lam}{\lambda}

\newcommand{\DD}{{\mathbb D}}

\newcommand{\CC}{{\mathbb C}}

\newcommand{\tr}{{\rm Tr}\,}
\newcommand{\ran}{{\rm ran}\,}
\newcommand{\ind}{\int_{\DD}}

\newcommand{\calL}{{\mathcal L}}

\begin{document}

\title{Frames and operators in Schatten classes}%

\date{\today}%

\author{Hu Bingyang, Le Hai Khoi, and Kehe Zhu}%

\address{(Hu \& Khoi) Division of Mathematical Sciences, School of Physical and 
Mathematical Sciences, Nanyang Technological University (NTU), 637371 Singapore }%
\email{BHU2@e.ntu.edu.sg; lhkhoi@ntu.edu.sg}

\address{(Zhu) Department of Mathematics, State University of New York at Albany, Albany, NY 12222, USA}%
\email{kzhu@math.albany.edu}


\begin{abstract}
Let $T$ be a compact operator on a separable Hilbert space $H$. We show that, for
$2\le p<\infty$, $T$ belongs to the Schatten class $S_p$ if and only if $\{\|Tf_n\|\}\in
\ell^p$ for \emph{every} frame $\{f_n\}$ in $H$; and for $0<p\le2$, $T$ belongs to $S_p$
if and only if $\{\|Tf_n\|\}\in\ell^p$ for \emph{some} frame $\{f_n\}$ in $H$. Similar
conditions are also obtained in terms of the sequence $\{\langle Tf_n,f_n\rangle\}$ and
the double-indexed sequence $\{\langle Tf_n,f_m\rangle\}$.
\end{abstract}

\maketitle

\section{Introduction}

Let $H$ be a separable Hilbert space. A sequence $\{f_n\}\subset H$ is called a \emph{frame} 
for $H$ if there exist positive constants $C_1$ and $C_2$ such that
$$C_1\|f\|^2 \le \sum_{n=1}^{\infty} |\langle f,f_n\rangle |^2 \le C_2\|f\|^2$$
for all $f\in H$. The numbers $C_1$ and $C_2$ are certainly not unique. The optimal upper 
constant, $\inf C_2$, will be called the upper frame bound for $\{f_n\}$. Similarly, the 
optimal lower constant, $\sup C_1$, will be called the lower frame bound for $\{f_n\}$. 
A frame is called \emph{tight} if its lower and upper frame bounds are the same.
Also, a frame is called \emph{Parseval} or \emph{normalized tight} if its lower and upper 
frame bounds are both $1$. See \cite{Chr} for an introduction to the theory of frames.

The \emph{singular values} or \emph{s-numbers} of a compact operator $T$ on $H$ are the 
square roots of the positive eigenvalues of the operator $T^*T$, where $T^*$ 
denotes the adjoint of $T$. Equivalently, this is the sequence of positive eigenvalues 
of $|T|=(T^*T)^{1/2}$. We always arrange the singular values of $T$, $\{\lambda_n\}$, 
such that $\lambda_1\ge\lambda_2\ge\lambda_3\ge\cdots$, with each eigenvalue
of multiplicity $k$ repeated $k$ times in the sequence.

Given $0<p<\infty$, the \emph{Schatten $p$-class} of $H$, denoted $S_p(H)$ or simply $S_p$, 
is defined as the space of all compact operators $T$ on $H$ with its singular value sequence $\{\lam_n\}$ belonging to $\ell^p$. It is well known that $S_p$ is a two-sided ideal in the
full algebra $\calL(H)$ of all bounded linear operators on $H$. Also, when equiped with
$$\|T\|_p=\left[\sum_{n=1}^\infty\lambda_n^p\right]^{\frac1p},$$
$S_p$ is a Banach space when $1\le p<\infty$ and a complete metric space when $0<p<1$.
Two special cases are especially interesting in operator theory: $S_1$ is called the
trace class and $S_2$ is the Hilbert-Schmidt class. See \cite{GK, Simon, Zhu} for 
basic properties of Schatten classes. 

Operators in Schatten classes can often be described by their action on orthonormal bases.
For example, a positive operator $T\in\calL(H)$ belongs to the trace class $S_1$ if and
only if $\sum\langle Te_n,e_n\rangle<\infty$, where $\{e_n\}$ is any given orthonormal basis 
for $H$. Similarly, an operator $T\in\calL(H)$ belongs to the Hilbert-Schmidt class $S_2$ if 
and only if $\sum\|Te_n\|^2<\infty$, where $\{e_n\}$ is any given orthonormal basis for $H$. 
See \cite{GK, Simon, Zhu} again for these and other related results.

It is clear that any orthonormal basis is a frame, with frame bounds equal to $1$.
The purpose of this article is to study Schatten class operators in terms of frames. 
We state our main results as follows.

\begin{thma}
Suppose $T$ is a compact operator on $H$ and $2\le p<\infty$. Then the following
conditions are equivalent.
\begin{enumerate}
\item[(a)] $T\in S_p$.
\item[(b)] $\{\|Te_n\|\}\in\ell^p$ for \emph{every} orthonormal basis $\{e_n\}$ in $H$.
\item[(c)] $\{\|Tf_n\|\}\in\ell^p$ for \emph{every} frame $\{f_n\}$ in $H$.
\end{enumerate}
Furthermore, we always have
$$\|T\|_p^p=\sup\sum_{n=1}^\infty\|Te_n\|^p=\sup\sum_{n=1}^\infty\|Tf_n\|^p,$$
where the first supremum is taken over all orthonormal bases $\{e_n\}$ and the second
supremum is taken over all frames $\{f_n\}$ with upper frame bound less than or equal to $1$.
\end{thma}

\begin{thmb}
Suppose $T$ is a bounded operator on $H$ and $0<p\le2$. Then the following
conditions are equivalent.
\begin{enumerate}
\item[(a)] $T\in S_p$.
\item[(b)] $\{\|Te_n\|\}\in\ell^p$ for \emph{some} orthonormal basis $\{e_n\}$ in $H$.
\item[(c)] $\{\|Tf_n\|\}\in\ell^p$ for \emph{some} frame $\{f_n\}$ in $H$.
\end{enumerate}
Furthermore, we always have
$$\|T\|_p^p=\inf\sum_{n=1}^\infty\|Te_n\|^p=\inf\sum_{n=1}^\infty\|f_n\|^{2-p}
\|Tf_n\|^p=\inf\sum_{n=1}^\infty\|Tf_n\|^p,$$
where the first infimum is taken over all orthonormal bases, the second infimum
is taken over all frames with lower frame bound greater than or equal to $1$, and the 
third infimum is taken over all Parseval frames $\{f_n\}$.
\end{thmb}

The conditions above concerning orthonormal basis are more or less well known to 
experts in the field. But the necessary \emph{and} sufficient conditions stated here do not 
seem to have appeared anywhere before. Partial statements in terms of orthonormal basis can 
be found in \cite{GK, Simon, Zhu} for example. We will include the treatment for orthonormal 
bases in the paper for the sake of completeness.

It is interesting to observe the sharp contrast between the cases $p\ge2$ and $p\le2$:
in the first case $\sup$ is used to compute the norm $\|T\|_p$, while in the second case
$\inf$ must be used. We will also construct examples to show that the cut-off at $p=2$ 
is necessary, and the result at the cut-off value $p=2$ is particularly nice. 

In addition to Theorems A and B, we will obtain corresponding results in terms of the 
sequence $\{\langle Tf_n,f_n\rangle\}$ and the double-indexed sequence 
$\{\langle Tf_n,f_k\rangle\}$. But in these cases it is sometimes \emph{necessary} to 
require additional assumptions on the operator $T$, such as $T$ being positive or
self-adjoint.

The relationship between frames and operators in Schatten classes has been studied by
several authors in the past few years. See \cite{ACP,Bal,Koo} and references therein. 
There is some overlap between the present paper and the papers just referenced. However,
the approach here is different, the results here are \emph{complete}, and the proofs here 
are simpler and more natural.

\section{The case when $p$ is large}

The description of operators in the Schatten class $S_p$ depends on the range of $p$.
In this section we focus on the case when $p$ large. We begin with the following lemma
which is well known to experts. This is the only result from the theory of frames that
we will use in the paper, so we include a short proof here for the reader's easy reference.

\begin{lem}
Suppose $\{e_n\}$ is an othonormal basis and $\{f_n\}$ is a frame for $H$. Then the operator
$A:H\to H$ defined by
$$A\left(\sum_{k=1}^\infty c_ke_k\right)=\sum_{k=1}^\infty c_kf_k$$
is a well-defined bounded linear operator. Furthermore, $\|A\|^2$ is between the lower and
upper frame bounds of $\{f_n\}$, and $AA^*$ is invertible on $H$.
\label{1}
\end{lem}

\begin{proof}
If $f$ is any vector in $H$, then
\begin{eqnarray*}
\left|\left\langle A\left(\sum_{k=1}^N c_ke_k\right),f\right\rangle\right|^2%
&=& \left|\sum_{k=1}^N c_k \langle f_k ,f\rangle\right|^2 \\
&\le& \sum_{k=1}^N |c_k|^2 \sum_{k=1}^N | \langle f_k, f \rangle |^2 \\
&\le& C_2 \|f\|^2\left\|\sum_{k=1}^N c_ke_k \right\|^2,
\end{eqnarray*}
where $C_2$ is the upper frame bound for $\{f_n\}$.

Therefore, by the Hahn-Banach theorem, $A$ extends to a bounded linear operator on $H$ 
with $\|A\|^2\le C_2$, namely,
$$A \left(\sum_{k=1}^{\infty} c_ke_k \right) = \sum_{k=1}^{\infty} c_kf_k,$$
where $\{c_k\}\in\ell^2$. If $\langle f,f_k\rangle=0$ for all $k$, then it follows from
the definition of frame that $f=0$. Therefore, $A$ has dense range.

For any vector $f\in H$, we have
$$A^*f=\sum_{n=1}^\infty\langle A^*f,e_n\rangle\,e_n=
\sum_{n=1}^\infty\langle f,Ae_n\rangle\,e_n=\sum_{n=1}^\infty\langle f,f_n\rangle\,e_n.$$
It follows that
$$\|A^*f\|^2=\sum_{n=1}^\infty|\langle f,f_n\rangle|^2\ge C_1\|f\|^2$$
for $f\in H$, where $C_1$ is the lower frame bound for $\{f_n\}$. This shows that
$\|A\|^2=\|A^*\|^2\ge C_1$ and $A^*$ is one-to-one and has closed range. Furthermore, 
for any $f\in H$, we have
$$C_1\|f\|^2\le\|A^*f\|^2=\langle AA^*f,f\rangle\le\|AA^*f\|\|f\|.$$
It follows that $\|AA^*f\|\ge C_1\|f\|$ for all $f\in H$, so that $AA^*$ is one-to-one
and has closed range. Since $\ran(AA^*)^\perp=\ker(AA^*)=(0)$, $AA^*$ must be onto. 
Therefore, $AA^*$ is invertible.
\end{proof}

As a consequence of the invertibility of $AA^*$, we see that the operator $A$ above is
actually onto. Therefore, every vector $f\in H$ admits a representation of the form
$$f=\sum_{n=1}^\infty c_nf_n,$$
where $\{c_n\}\in\ell^2$. Note that $A$ is generally not one-to-one. For example, a frame 
may contain a certain vector that is repeated a finite number of times. In this case, the
associated operator $A$ is obviously not one-to-one.

\begin{thm}
Suppose $T$ is a compact operator on $H$ and $2\le p<\infty$. Then
the following statements are equivalent.
\begin{enumerate}
\item[(a)] $T$ is in the Schatten class $S_p$.
\item[(b)] $\{\|Te_n\|\}\in\ell^p$ for \emph{every} orthonormal basis $\{e_n\}$ in $H$.
\item[(c)] $\{\|Tf_n\|\}\in\ell^p$ for \emph{every} frame $\{f_n\}$ in $H$.
\end{enumerate}
Moreover, we always have
$$\|T\|^p_p=\sup\sum_{n=1}^\infty\|Te_n\|^p=\sup\sum_{n=1}^\infty\|Tf_n\|^p,$$
where the first supremum is taken over all orthonormal bases $\{e_n\}$ and the second
supremum is taken over all frames $\{f_n\}$ with upper frame bound less than or
equal to $1$.
\label{2}
\end{thm}

\begin{proof}
The equivalence of conditions (a) and (b) is well known. See Theorem 1.33 of \cite{Zhu} for 
example. Note that Theorem 1.33 of \cite{Zhu} was stated and proved in terms of orthonormal 
sets. Since every orthonormal set can be expanded to an orthonormal basis, the result 
remains true when the phrase ``orthonormal sets'' is replaced by ``orthonormal bases''.

Since every orthonormal basis is a frame, it is trivial that condition (c) implies (b).

To prove that (a) and (b) together imply (c), we fix an orthonormal basis $\{e_n\}$ and
a frame $\{f_n\}$ for $H$ and consider the operator $A$ defined in Lemma~\ref{1}. If
$T$ is in $S_p$, then so is the operator $S=TA$. Apply condition (b) to the operator $S$,
we obtain
$$\sum_{n=1}^\infty\|Tf_n\|^p=\sum_{n=1}^\infty\|TAe_n\|^p
=\sum_{n=1}^\infty\|Se_n\|^p<\infty.$$
This completes the proof of the equivalence of conditions (a), (b), and (c).

The equality 
$$\|T\|_p^p=\sup\sum_{n=1}^\infty\|Te_n\|^p$$
was established in Theorem 1.33 of \cite{Zhu}. Since every orthonormal basis is a frame
with frame bounds equal to $1$, we clearly have
$$\|T\|_p^p=\sup\sum_{n=1}^\infty\|Te_n\|^p\le\sup\sum_{n=1}^\infty\|Tf_n\|^p.$$
This along with the arguments in the previous paragraph shows that
$$\sum_{n=1}^\infty\|Tf_n\|^p=\sum_{n=1}^\infty\|TAe_n\|^p\le\|TA\|_p^p.$$
It is well known (see \cite{GK,Zhu} for example) that $\|TA\|_p\le\|T\|_p\|A\|$, which 
combined with the estimate for $\|A\|$ in Lemma~\ref{1} shows that $\|TA\|_p\le\|T\|_p$
whenever $\{f_n\}$ has upper frame bound less than or equal to $1$. This shows that
$$\sup\sum_{n=1}^\infty\|Tf_n\|^p_p\le\|T\|_p^p,$$
and completes the proof of the theorem.
\end{proof}

It is possible to obtain a version of Theorem~\ref{2} without the a priori assumption that
$T$ be compact. This will be done using an approximation argument based on the following lemma.

\begin{lem}
Suppose $T$ and $T_k$, $k\ge1$, are bounded linear operators on $H$. If $1<p<\infty$,
$T_k\to T$ in the weak operator topology, and $\|T_k\|_p\le C$ for some constant $C$ 
and all $k\ge1$, then $\|T\|_p\le C$.
\label{3}
\end{lem}

\begin{proof}
Let $S$ be a finite-rank operator and $\{e_n\}$ be an orthonormal basis of $H$ such that
$TS(e_n)=0$ for all but a finite number of $n$. Then
$$\tr(TS)=\sum_{n=1}^\infty\langle TSe_n,e_n\rangle=\lim_{k\to\infty}
\sum_{n=1}^\infty\langle T_kSe_n,e_n\rangle=\lim_{k\to\infty}\tr(T_kS).$$
Since the Banach dual of $S_p$ is $S_q$, $1/p+1/q=1$, under the pairing induced by 
the trace, we have
$$|\tr(T_kS)|\le\|T_k\|_p\|S\|_q\le C\|S\|_q,\qquad k\ge1.$$
It follows that $|\tr(TS)|\le C\|S\|_q$ for all finite rank 
operators $S$. Since the set of finite rank operators is dense in $S_q$, we conclude
that $T\in S_p$ and $\|T\|_p\le C$.
\end{proof}

\begin{thm}
When $2\le p<\infty$, the following conditions are equivalent for any bounded linear 
operator $T$ on $H$:
\begin{enumerate}
\item[(i)] $T\in S_p$.
\item[(ii)] There exists a constant $C>0$ such that
$$\sum_{n=1}^\infty\|Te_n\|^p\le C$$
for \emph{every} orthonormal basis $\{e_n\}$.
\item[(iii)] There exists a constant $C>0$ such that
$$\sum_{n=1}^\infty\|Tf_n\|^p\le C$$
for \emph{every} frame $\{f_n\}$ with upper frame bound no greater than $1$.
\end{enumerate}
\label{4}
\end{thm}

\begin{proof}
All we have to show here is that condition (ii) implies (i) without the a priori assumption
that $T$ be compact. This can be done with the help of an approximation argument. More
specifically, we fix an increasing sequence $\{P_k\}$ of finite-rank projections such that
$\{P_k\}$ converges to the identity operator in the strong operator topology and let
$T_k=P_kT$ for $k\ge1$. Each $T_k$ is a finite-rank operator, so condition (ii) along
with Theorem~\ref{2} gives
$$\|T_k\|^p_p=\sup\sum_{n=1}^\infty\|T_ke_n\|^p\le\sup\sum_{n=1}^\infty\|Te_n\|^p\le C$$
for all $k\ge1$, where the suprema are taken over all orthonormal bases $\{e_n\}$. Since 
$T_k\to T$ in the strong operator topology, it follows from Lemma~\ref{3} that 
$\|T\|_p^p\le C$.
\end{proof}

It is natural to ask the following question: suppose $p\ge2$ and $\{\|Tf_n\|\}\in\ell^p$ for \emph{some} frame $\{f_n\}$, does it imply that $T\in S_p$? The answer is yes for $p=2$ but
no for $p>2$. We will get back to the case $p=2$ in Section 4 but will now settle the case
$p>2$. 

\begin{prop}
Suppose $2<p<\infty$, $\varepsilon>0$ (not necessarily small), and $T$ is any operator 
in $S_{p+\varepsilon}-S_p$. Then there exists a frame $\{f_n\}$ for $H$ such that 
$\{\|Tf_n\|\}\in\ell^p$.
\label{5}
\end{prop}

\begin{proof}
Suppose that
$$Tx=\sum_{n=1}^\infty\lambda_n\langle x,e_n\rangle\,\sigma_n$$
is the canonical decomposition of $T$, where $\{\lambda_n\}$ is the singular value sequece 
of $T$ which is arranged in nonincreasing order and repeated according to multiplicity. 
Thus we have $\{\lambda_n\}\in\ell^{p+\varepsilon}$ but $\{\lambda_n\}\not\in\ell^p$.

Let $\{e_n'\}$ denote an orthonormal basis for $\ker(T)$. Then $\{e_n\}\cup\{e_n'\}$ is an
orthonormal basis for $H$. In fact, for any vector $x\in H$, we have $\langle x,e_n\rangle=0$ 
for every $n$ if and only if $Tx=0$. Therefore, $\{e_n\}^\perp=\{e_n'\}$.

For every $n\ge1$ choose a positive number $\delta_n$ such that $\delta_n^{p-2}
=\lambda_n^\varepsilon$. Since $p>2$, we have $\delta_n\to0$ as $n\to\infty$, so we
can choose a sequence $\{N_n\}$ of positive integers such that $N_n\delta_n^2\sim1$ 
as $n\to\infty$. In other words, there exist positive constants $c$
and $C$ such that $c\le N_n\delta_n^2\le C$ for all $n\ge1$. Let $\{f_n\}$ be the sequence 
consisting of all vectors in $\{e_n'\}$, plus $N_1$ copies of the vector $\delta_1e_1$, 
plus $N_2$ copies of the vector $\delta_2e_2$, and so on.

For any vector $f\in H$, we have
\begin{eqnarray*}
\sum_{n=1}^\infty|\langle f,f_n\rangle|^2&=&\sum_{n=1}^\infty N_n|\langle f,
\delta_ne_n\rangle|^2+\sum_{n=1}^\infty|\langle f,e_n'\rangle|^2\\
&=&\sum_{n=1}^\infty N_n\delta_n^2|\langle f,e_n\rangle|^2+\sum_{n=1}^\infty
|\langle f,e_n'\rangle|^2\\
&\sim&\sum_{n=1}^\infty|\langle f,e_n\rangle|^2+\sum_{n=1}^\infty|\langle f,e_n'\rangle|^2\\
&=&\|f\|^2.
\end{eqnarray*}
This shows that $\{f_n\}$ is a frame for $H$.

On the other hand,
\begin{eqnarray*}
\sum_{n=1}^\infty\|Tf_n\|^p&=&\sum_{n=1}^\infty N_n\|T(\delta_ne_n)\|^p
=\sum_{n=1}^\infty N_n\delta_n^p\|Te_n\|^p\\
&=&\sum_{n=1}^\infty N_n\delta_n^p\lambda_n^p
=\sum_{n=1}^\infty N_n\delta_n^2\delta_n^{p-2}\lambda_n^p\\
&\sim&\sum_{n=1}^\infty\delta_n^{p-2}\lambda_n^p
=\sum_{n=1}^\infty\lambda_n^{p+\varepsilon}<\infty.
\end{eqnarray*}
This completes the proof of the proposition.
\end{proof}

We also derive a version of the above proposition in terms of orthonormal bases.

\begin{prop}
Suppose $2<p<\infty$, $\varepsilon>0$, and $\{e_n\}$ is any orthonormal basis for $H$.
Then there exists an operator $S\in S_{p+\varepsilon}-S_p$ such that $\{\|Se_n\|\}\in\ell^p$.
\label{6}
\end{prop}

\begin{proof}
Fix any operator $T\in S_{p+\varepsilon}-S_p$ and use Proposition~\ref{5} to select a 
frame $\{f_n\}$ such that $\{\|Tf_n\|\}\in\ell^p$. Let $A$ be the operator on $H$ defined
by $Ae_n=f_n$, $n\ge1$. By Lemma~\ref{1}, the operator $A$ is bounded and the operator
$AA^*$ is invertible. Let $S=TA$. Then $S\in S_{p+\varepsilon}$ because $S_{p+\varepsilon}$
is a two-sided ideal in the full operator algebra $\calL(H)$. Since $AA^*$ is invertible,
we have $S\not\in S_p$ as well. Otherwise, the operator $T(AA^*)=SA^*$ would be in
$S_p$. Multiplying from the right by $(AA^*)^{-1}$ and using the fact that
$S_p$ is a two-sided ideal in $\calL(H)$ again, we would then obtain that $T$ is in $S_p$,
a contradiction. Therefore, $S\in S_{p+\varepsilon}-S_p$ and
$$\sum_{n=1}^\infty\|Se_n\|^p=\sum_{n=1}^\infty\|TAe_n\|^p=\sum_{n=1}^\infty
\|Tf_n\|^p<\infty.$$
This completes the proof of the proposition.
\end{proof}

Characterizations of Schatten classes can also be given in terms of the
sequence $\{\langle Tf_n,f_n\rangle\}$ and the double-indexed sequence $\{\langle Tf_n,
f_k\rangle\}$. We now proceed to the characterization of $S_p$ based on the sequence 
$\{\langle Tf_n,f_n\rangle\}$.

\begin{thm}
Suppose $1\le p<\infty$ and $S$ is a compact operator on $H$. Then the following
conditions are equivalent.
\begin{enumerate}
\item[(a)] $S$ belongs to $S_p$.
\item[(b)] $\{\langle Se_n,e_n\rangle\}\in\ell^p$ for \emph{every} orthonormal 
basis $\{e_n\}$ in $H$.
\item[(c)] $\{\langle Sf_n,f_n\rangle\}\in\ell^p$ for \emph{every} frame $\{f_n\}$ in $H$. 
\end{enumerate}
Furthermore, if $S$ is self-adjoint, then
$$\|S\|_p^p=\sup\sum_{n=1}^\infty|\langle Se_n,e_n\rangle|^p=\sup\sum_{n=1}^\infty
|\langle Sf_n,f_n\rangle|^p,$$
where the first supremum is taken over all orthonormal bases and the second supremum is
taken over all frames with upper frame bound less than or equal to $1$.
\label{7}
\end{thm}

\begin{proof}
The equivalence of (a) and (b) follows from Theorem 1.27 in \cite{Zhu}. Note again that
Theorem 1.27 in \cite{Zhu} is stated in terms of orthonormal sets. Since every orthonormal
set can be expanded to an orthonormal basis, we see that Theorem 1.27 in \cite{Zhu} remains
valid when the phrase``orthonormal sets'' is replaced by ``orthonormal bases''.

It is trivial that (c) implies (b).

To prove that (a) implies (c), first assume that $S$ is positive. In this case, we can
write $S=T^*T$, where $T$ is the square root of $S$. Then $\langle Sf_n,f_n\rangle
=\|Tf_n\|^2$ and the desired result follows from Theorem~\ref{2} and the fact that 
$S\in S_p$ if and only if $T\in S_{2p}$. When $S$ is not necessarily positive, it is 
well known that we can write 
$$S=(S_1-S_2)+i(S_3-S_4),$$
where each $S_k$ is positive and belongs to $S_p$. By the already proved case for positive
operators, $\{\langle S_kf_n,f_n\rangle\}\in\ell^p$ for each $1\le k\le4$. It follows that
$\{\langle Sf_n,f_n\rangle\}\in\ell^p$.

It follows from the canonical decomposition for self-adjoint compact operators and the
fact that every orthonormal basis is a frame with frame bounds $1$ that
$$\|S\|_p^p\le\sup\sum_{n=1}^\infty|\langle Se_n,e_n\rangle|^p
\le\sup\sum_{n=1}^\infty|\langle Sf_n,f_n\rangle|^p.$$
If $\{f_n\}$ is a frame with upper frame bound $1$, then by the norm estimate for
$A$ in Lemma~\ref{1} we have $\|f_n\|\le1$ for every $n$, which together with 
Theorem~\ref{2} gives
$$\sum_{n=1}^\infty|\langle Sf_n,f_n\rangle|^p\le\sum_{n=1}^\infty
\|Sf_n\|^p\le\|S\|_p^p.$$
This shows that
$$\sup\sum_{n=1}^\infty|\langle Sf_n,f_n\rangle|^p\le\|S\|_p^p,$$
and completes the proof of the theorem.
\end{proof}

Note that the second assertion in Theorem~\ref{7} concerning the norm of $S$ in $S_p$ is
false for operators that are not necessarily self-adjoint. A counter-example can be found 
on page 22 of \cite{Zhu}. Nevertheless, using the fact that every operator $T$ admits a
canonical decomposition $T=T_1+iT_2$ with
$$|\langle Tf,f\rangle|^2=|\langle T_1f,f\rangle^2+|\langle T_2f,f\rangle|^2,$$
where
$$T_1=\frac{T+T^*}2,\qquad T_2=\frac{T-T^*}{2i},$$
are self-adjoint, we easily show that there still exists a positive constant $C$ such that
$$C^{-1}\|T\|_p^p\le\sup\sum_{n=1}^\infty|\langle Te_n,e_n\rangle|^p\le C\|T\|_p^p$$
for all operators $T\in S_p$, where the supremum is taken over all orthonormal bases
$\{e_n\}$. See the second part of the proof of Theorem 1.27 in \cite{Zhu}. Similarly, we have
$$C^{-1}\|T\|_p^p\le\sup\sum_{n=1}^\infty|\langle Tf_n,f_n\rangle|^p\le C\|T\|_p^p$$
for all operators $T\in S_p$, where the supremum is taken over all frames $\{f_n\}$ with
upper frame bound less than or equal to $1$.

If we remove the a priori assumption that $T$ be compact, we obtain the following slightly
different version of Theorem~\ref{7}.

\begin{thm}
If $1\le p<\infty$ and $S$ is a bounded linear operator on $H$, then the following 
conditions are equivalent:
\begin{enumerate}
\item[(i)] $S\in S_p$.
\item[(ii)] There exists a positive constant $C$ such that
$$\sum_{n=1}^\infty|\langle Se_n,e_n\rangle|^p\le C$$
for \emph{every} orthonormal basis $\{e_n\}$.
\item[(iii)] There exists a positive constant $C$ such that
$$\sum_{n=1}^\infty|\langle Sf_n,f_n\rangle|^p\le C$$
for \emph{every} frame $\{f_n\}$ with upper frame bound less than or equal to $1$.
\end{enumerate}
\label{8}
\end{thm}

\begin{proof}
This follows from Theorem~\ref{7}, the remarks immediately following Theorem~\ref{7}, 
and the same approximation argument used in the proof of Theorem~\ref{4}.
\end{proof}

\begin{prop}
If $1<p<\infty$, $\varepsilon>0$, and $S\in S_{p+\varepsilon}-S_p$ is positive, then 
there exists \emph{some} frame $\{f_n\}$ such that $\{\langle Sf_n,f_n\rangle\}\in\ell^p$.
\label{9}
\end{prop}

\begin{proof}
Write $S=T^*T$, where $T=\sqrt S$. Then $T\in S_{2p+2\varepsilon}-S_{2p}$ and 
$\|Tf_n\|^{2p}=\langle Sf_n,f_n\rangle^p$. The desired result then follows from 
Proposition~\ref{5}.
\end{proof}

\begin{prop}
If $1<p<\infty$, $\varepsilon>0$, and $\{e_n\}$ is an orthonormal basis for $H$, then
there exists a positive operator $S\in S_{p+\varepsilon}-S_p$ such that 
$\{\langle Se_n,e_n\rangle\}\in\ell^p$.
\label{10}
\end{prop}

\begin{proof}
By Proposition \ref{6}, there exists an operator $T\in S_{2p+2\varepsilon}-S_{2p}$
such that $\{\|Te_n\|\}\in\ell^{2p}$. Let $S=T^*T$. Then $S\in S_{p+\varepsilon}-S_p$
and the sequence $\langle Se_n,e_n\rangle=\|Te_n\|^2$ belongs to $\ell^p$.
\end{proof}

Next we proceed to the characterization of operators in Schatten classes in terms of
the double-indexed sequence $\{\langle Tf_n,f_k\rangle\}$. We need the following lemma.

\begin{lem}
For any frame $\{f_n\}$ in $H$ there exist positive constants $c$ and $C$ with the 
following properties.
\begin{enumerate}
\item[(a)] If $2\le p<\infty$, then
$$\sum_{n=1}^\infty\sum_{k=1}^\infty|\langle Tf_n,f_k\rangle|^p\le C\sum_{n=1}^\infty
\|Tf_n\|^p$$
for all bounded linear operators $T$ on $H$.
\item[(b)] If $0<p\le2$, then
$$\sum_{n=1}^\infty\sum_{k=1}^\infty|\langle Tf_n,f_k\rangle|^p\ge c\sum_{n=1}^\infty
\|Tf_n\|^p$$
for all bounded linear operators on $H$.
\end{enumerate}
\label{11}
\end{lem}

\begin{proof}
The desired estimates follow from H\"older's inequality and the definition of frames.

For $2\le p<\infty$, we have $0<2/p\le1$, so
$$\left[\sum_{k=1}^\infty|\langle Tf_n,f_k\rangle|^p\right]^{\frac2p}
\le\sum_{k=1}^\infty|\langle Tf_n,f_k\rangle|^2\le C_1\|Tf_n\|^2,$$
where $C_1$ is upper frame bound for $\{f_n\}$. It follows that for $C=C_1^{p/2}$ we have
$$\sum_{k=1}^\infty|\langle Tf_n,f_k\rangle|^p\le C\|Tf_n\|^p,\qquad n\ge1,$$
so that
$$\sum_{n=1}^\infty\sum_{k=1}^\infty|\langle Tf_n,f_k\rangle|^p
\le C\sum_{n=1}^\infty\|Tf_n\|^p.$$

Similarly, for $0<p\le 2$, we have $0<p/2\le1$, so
\begin{eqnarray*}
\sum_{n=1}^\infty\sum_{k=1}^\infty|\langle Tf_n,f_k\rangle|^p&=&
\sum_{n=1}^\infty\sum_{k=1}^\infty\left(|\langle Tf_n,f_k\rangle|^2\right)^{\frac p2}\\
&\ge&\sum_{n=1}^\infty\left[\sum_{k=1}^\infty|\langle Tf_n,f_k\rangle|^2\right]^{\frac p2}\\
&\ge&C_2^{p/2}\sum_{n=1}^\infty\left[\|Tf_n\|^2\right]^{\frac p2}\\
&=&c\sum_{n=1}^\infty\|Tf_n\|^p,
\end{eqnarray*}
where $C_2$ is the lower frame bound for $\{f_n\}$ and $c=C_2^{p/2}$.
\end{proof}

It is clear that if $\{f_n\}$ happens to be an orthonormal basis, then both $C$ and $c$ can
be taken to be $1$ in Lemma~\ref{11}.

\begin{thm}
Suppose $T$ is a compact operator on $H$ and $2\le p<\infty$. Then the following
conditions are equivalent.
\begin{enumerate}
\item[(a)] $T\in S_p$.
\item[(b)] The condition
$$\sum_{n=1}^\infty\sum_{k=1}^\infty|\langle Te_n,e_k\rangle|^p<\infty$$
for \emph{every} orthonormal basis $\{e_n\}$ in $H$.
\item[(c)] The condition
$$\sum_{n=1}^\infty\sum_{k=1}^\infty|\langle Tf_n,f_k\rangle|^p<\infty$$
holds for \emph{every} frame $\{f_n\}$ in $H$.
\end{enumerate}
Furthermore, there exists a positive constant $c$ such that
$$c\|T\|_p^p\le\sup\sum_{n=1}^\infty\sum_{k=1}^\infty|\langle Te_n,e_k\rangle|^p
\le\sup\sum_{n=1}^\infty\sum_{k=1}^\infty|\langle Tf_n,f_k\rangle|^p\le\|T\|_p^p,$$
where the first supremum is taken over all orthonormal bases $\{e_n\}$ and the second
supremum is taken over all frames $\{f_n\}$ with upper frame bound less than or equal
to $1$.
\label{12}
\end{thm}

\begin{proof}
That condition (a) implies (c) follows from Theorem~\ref{2} and part (a) 
of Lemma~\ref{11}. Since every orthonormal basis is a frame, it is trivial that 
condition (c) implies (b). 

It remains to show that condition (b) implies (a). So we assume that condition (b) holds for
an operator $T$. It is clear that condition (b) also holds for $T^*$, which implies that
condition (b) holds for $T+T^*$ and $T-T^*$ as well. Write
$T=T_1+iT_2$, where $T_1=(T+T^*)/2$ and $T_2=(T-T^*)/(2i)$, and apply condition (b) to the
self-adjoint operators $T_1$ and $T_2$, we may as well assume that $T$ is already self-adjoint.

But if $T$ is self-adjoint, its canonical decomposition takes the form
$$Tx=\sum_{n=1}^\infty\lambda_n\langle x,e_n\rangle\,e_n,$$
where $\{\lambda_n\}$ is the singular value sequence of $T$ and $\{e_n\}$ is an orthonormal
set. If $\{\sigma_n\}$ is an orthonormal basis for $\ker(T)$, then $\{e_n'\}=\{e_n\}\cup
\{\sigma_n\}$ is an orthonormal basis for $H$. Therefore, it follows from condition (b) 
and the relation $\ran(T)^\perp=\ran(T^*)^\perp=\ker(T)$ that
$$\sum_{n=1}^\infty\lambda_n^p=\sum_{n=1}^\infty\sum_{k=1}^\infty|\langle Te_n,e_k\rangle|^p
=\sum_{n=1}^\infty\sum_{k=1}^\infty|\langle Te_n',e_k'\rangle|^p<\infty.$$
The first norm estimate follows from the decomposition $T=T_1+iT_2$ of $T$ into a linear
combination of self-adjoint operators and the canonical decomposition of self-adjoint
compact operators. The second norm estimate is trivial. The third norm estimate follows 
from Theorem~\ref{2} and part (a) of Lemma~\ref{11}.
\end{proof}

Note that if $T$ is self-adjoint, then the proof above actually shows
$$\|T\|_p^p=\sup\sum_{n=1}^\infty\sum_{k=1}^\infty|\langle Te_n,e_k\rangle|^p
=\sup\sum_{n=1}^\infty\sum_{k=1}^\infty|\langle Tf_n,f_k\rangle|^p.$$
We are not sure if this holds for general operators as well.

Once again, if we do not make the a priori assumption that $T$ be compact, then 
Theorem~\ref{12} should be modified as follows.

\begin{thm}
For $2\le p<\infty$ and any bounded linear operator $T$ the following conditions 
are equivalent:
\begin{enumerate}
\item[(i)] $T\in S_p$.
\item[(ii)] There exists a positive constant $C$ such that
$$\sum_{n=1}^\infty\sum_{k=1}^\infty|\langle Te_n,e_k\rangle|^p\le C$$
for \emph{every} orthonormal basis $\{e_n\}$.
\item[(iii)] There exists a positive constant $C$ such that
$$\sum_{n=1}^\infty\sum_{k=1}^\infty|\langle Tf_n,f_k\rangle|^p\le C$$
for \emph{every} frame $\{f_n\}$ with upper frame bound less than or equal to $1$.
\end{enumerate}
\label{13}
\end{thm}

\begin{proof}
This follows from Theorem~\ref{12} and the approximation argument used in the proof 
of Theorems~\ref{4} and \ref{8}.
\end{proof}

\begin{prop}
Let $2<p<\infty$, $\varepsilon>0$, and $T\in S_{p+\varepsilon}-S_p$. There exists a 
frame $\{f_n\}$ such that
$$\sum_{k=1}^\infty\sum_{n=1}^\infty|\langle Tf_n,f_k\rangle|^p<\infty.$$
\label{14}
\end{prop}

\begin{proof}
This follows from Proposition~\ref{5} and part (a) of Lemma~\ref{11}.
\end{proof}

\begin{prop}
Suppose $2<p<\infty$, $\varepsilon>0$, and $\{e_n\}$ is an orthonormal basis for $H$.
Then there exists an operator $T\in S_{p+\varepsilon}-S_p$ such that
$$\sum_{k=1}^\infty\sum_{n=1}^\infty|\langle Te_n,e_k\rangle|^p<\infty.$$
\label{15}
\end{prop}

\begin{proof}
This follows from Proposition~\ref{6} and part (a) of Lemma~\ref{11}.
\end{proof}

\section{The case when $p$ is small}

We begin this section with a simple example to show that the characterizations for 
operators in Schatten classes $S_p$ obtained in the previous section for $2\le p<\infty$
are not true for the range $0<p<2$.

Fix any orthonormal basis $\{e_n\}$ and consider the vector
$$h=\sum_{n=1}^\infty\frac{e_n}{\sqrt n\,\log(n+1)}$$
in $H$. Define a rank one operator $T$ on $H$ by $Tx=\langle x,h\rangle\,h$. We have
$$Te_n=\langle e_n,h\rangle\,h=\frac h{\sqrt n\,\log(n+1)},\qquad n\ge1.$$
It follows that
$$\sum_{n=1}^\infty\|Te_n\|^p=\|h\|^p\sum_{n=1}^\infty\frac1{[\sqrt n\,\log(n+1)]^p}
=\infty$$
for any $0<p<2$. This shows that the characterizations obtained in Theorem~\ref{2}
are no longer true for \emph{any} $0<p<2$.

Later in this section we will actually show that for \emph{any} operator $T\in S_p$,
$0<p<2$, there exists a frame $\{f_n\}$ such that $\{\|Tf_n\|\}\not\in\ell^p$. 
Nevertheless, there is still a nice characterization for operators in $S_p$, $0<p\le2$,
in terms of orthonormal bases and frames.

\begin{thm}
Suppose $T$ is a positive operator on $H$ and $0<p\le1$. Then the following 
conditions are equivalent.
\begin{enumerate}
\item[(a)] $T\in S_p$.
\item[(b)] $\{\langle Te_n,e_n\rangle\}\in\ell^p$ for \emph{some} orthonormal 
basis $\{e_n\}$ in $H$.
\item[(c)]$\{\langle Tf_n,f_n\rangle\}\in\ell^p$ for \emph{some} frame $\{f_n\}$ in $H$.
\end{enumerate}
Furthermore, we have
\begin{eqnarray*}
\|T\|^p_p&=&\inf\sum_{n=1}^\infty|\langle Te_n,e_n\rangle|^p\\
&=&\inf\sum_{n=1}^\infty\|f_n\|^{2(1-p)}\langle Tf_n,f_n\rangle^p\\
&=&\inf\sum_{n=1}^\infty\langle Tf_n,f_n\rangle^p,
\end{eqnarray*}
where the first infimum is taken over all orthonormal bases, the second
infimum is taken over all frames with lower frame bound greater than
or equal to $1$, and the third infimum is taken over all Parseval frames.
\label{16}
\end{thm}

\begin{proof}
If $T$ is positive and in $S_p$, then its canonical decomposition takes the form
$$Tx=\sum_{n=1}^\infty\lambda_n\langle x,\sigma_n\rangle\,\sigma_n,$$
where $\{\lambda_n\}\in\ell^p$ is the singular value sequence of $T$ and 
$\{\sigma_n\}$ is an orthonormal set in $H$. Since each $\lambda_n$ is positive, we 
have $Tx=0$ if and only if $\langle x,\sigma_n\rangle=0$ for every $n$. Therefore, 
$\ker(T)=\{\sigma_n\}^\perp$. If $\{\sigma_n'\}$ is an orthonormal basis for $\ker(T)$, 
then $\{e_n\}=\{\sigma_n\}\cup\{\sigma_n'\}$ is an orthonormal basis for $H$. Since 
$T(\sigma_n)=\lambda_n \sigma_n$ for every $n$ and $\{\lambda_n\}\in\ell^p$, we have
$$\sum_{n=1}^\infty|\langle Te_n,e_n\rangle|^p=\sum_{n=1}^\infty|\langle T\sigma_n,
\sigma_n\rangle|^p=\sum_{n=1}^\infty\lambda_n^p=\|T\|_p^p<\infty.$$
This shows that condition (a) implies (b).

To prove that condition (b) implies (a), we use Theorem 1.26 and part (b) of 
Proposition 1.31 in \cite{Zhu}. More specifically,
$$\|T\|_p^p=\|T^p\|_1^p=\sum_{n=1}^\infty\langle T^pe_n,e_n\rangle
\le\sum_{n=1}^\infty\langle Te_n,e_n\rangle^p<\infty$$
whenever $\{\langle Te_n,e_n\rangle\}\in\ell^p$. This shows that condition (b) implies (a), 
so conditions (a) and (b) are equivalent for any positive operator on $H$ and
$$\|T\|_p^p=\inf\sum_{n=1}^\infty\langle Te_n,e_n\rangle^p,$$
where the infimum is taken over all orthonormal bases $\{e_n\}$.

Since every orthonormal basis is a frame with both upper and lower frame bounds equal 
to $1$, it is trivial that condition (b) implies (c), and
$$\inf\sum_{n=1}^\infty\|f_n\|^{2(1-p)}\langle Tf_n,f_n\rangle^p\le\inf\sum_{n=1}^\infty
\langle Te_n,e_n\rangle^p,$$
where the first infimum is taken over all frames $\{f_n\}$ with lower frame bound at
least $1$ and the second infimum is taken over all orthonormal bases $\{e_n\}$.

Finally, we assume that $\{\langle Tf_n,f_n\rangle\}\in\ell^p$ for some frame $\{f_n\}$.
Fix any orthonormal basis $\{e_n\}$, let $A$ be the operator defined in Lemma~\ref{1}, and 
set $S=A^*TA$. Then $S$ is positive again and
$$\langle Se_n,e_n\rangle=\langle A^*TAe_n,e_n\rangle=\langle TAe_n,Ae_n\rangle
=\langle Tf_n,f_n\rangle.$$
Thus $\{\langle Se_n,e_n\rangle\}\in\ell^p$. By the equivalence of (a) and (b),
$S$ is in $S_p$. Since $S_p$ is a two-sided ideal in $\calL(H)$, the operator
$(AA^*)T(AA^*)=ASA^*$ also belongs to $S_p$. Multiply from both sides by $(AA^*)^{-1}$ (see
Lemma~\ref{1}), we conclude that $T$ is in $S_p$ as well. This shows that condition (c) 
implies (a), and completes the proof of the equivalence of (a), (b), and (c).

We now proceed to prove that
$$\|T\|_p^p\le\sum_{n=1}^\infty\|f_n\|^{2(1-p)}\langle Tf_n,f_n\rangle^p$$
whenever $\{f_n\}$ is a frame with lower frame bound greater than or equal to $1$. To 
this end, we fix such a frame and, without loss of generality, assume that the right-hand
side above is finite (otherwise, the desired inequality is trivial) and $f_n\not=0$
for each $n$. By part (b) of Proposition 1.31 in \cite{Zhu}, we have
\begin{eqnarray*}
+\infty&>&\sum_{n=1}^\infty\|f_n\|^{2(1-p)}\langle Tf_n,f_n\rangle^p\\
&=&\sum_{n=1}^\infty\|f_n\|^2\left\langle T\frac{f_n}{\|f_n\|},
\frac{f_n}{\|f_n\|}\right\rangle^p\\
&\ge&\sum_{n=1}^\infty\|f_n\|^2\left\langle T^p\frac{f_n}{\|f_n\|},
\frac{f_n}{\|f_n\|}\right\rangle\\
&=&\sum_{n=1}^\infty\langle T^pf_n,f_n\rangle.
\end{eqnarray*}
By the equivalence of (a) and (c), the operator $T^p$ is in the trace class. If
$$T^px=\sum_{k=1}^\infty\mu_k\langle x,\sigma_k\rangle\,\sigma_k$$
is the canonical decomposition for $T^p$, then
$$\langle T^pf_n,f_n\rangle=\sum_{k=1}^\infty\mu_k|\langle f_n,\sigma_k\rangle|^2$$
for every $n$. Since the lower frame bound of $\{f_n\}$ is greater than or equal to $1$,
it follows from Fubini's theorem and Theorem 1.26 in \cite{Zhu} that
\begin{eqnarray*}
\sum_{n=1}^\infty\langle T^pf_n,f_n\rangle^p&=&\sum_{k=1}^\infty\mu_k
\sum_{n=1}^\infty|\langle f_n,\sigma_k\rangle|^2\\
&\ge&\sum_{k=1}^\infty\mu_k\|\sigma_k\|^2=\sum_{k=1}^\infty\mu_k\\
&=&\|T^p\|_1=\|T\|_p^p.
\end{eqnarray*}
This completes the proof that
$$\|T\|_p^p=\inf\sum_{n=1}^\infty\langle Te_n,e_n\rangle^p=\inf\sum_{n=1}^\infty
\|f_n\|^{2(1-p)}\langle Tf_n,f_n\rangle^p,$$
where the first infimum is taken over all orthonormal bases and the second infimum is
taken over all frames with lower frame bound at least $1$.

Finally, if $\{f_n\}$ is a Parseval frame, then by the norm estimate for $A$ in Lemma~\ref{1},
we have $\|f_n\|\le1$ for every $n$. It follows that
$$\sum_{n=1}^\infty\|f_n\|^{2(1-p)}\langle Tf_n,f_n\rangle^p\le\sum_{n=1}^\infty
\langle Tf_n,f_n\rangle^p.$$
Combining this with the fact that every orthonormal basis is a Parseval frame, we obtain
$$\|T\|_p^p=\inf\sum_{n=1}^\infty\langle Tf_n,f_n\rangle^p,$$
where the infimum is taken over all Parseval frames $\{f_n\}$.
\end{proof}

Note that, unlike Theorem \ref{7}, we need to make the additional assumption that $T$ be
positive here. Without any extra assumption, Theorem~\ref{16} will be false. For
example, if $\{e_n\}$ is any fixed orthonormal basis for $H$ and $T$ is the unilateral shift
operator defined by $T(e_n)=e_{n+1}$, $n\ge1$. Then it is clear that 
$\{\langle Te_n,e_n\rangle\}\in\ell^p$ for any $p>0$, but $T$ is not even compact.

As a consequence of Theorem \ref{16} we obtain the following.

\begin{thm}
Suppose $T$ is a bounded linear operator on $H$ and $0<p\le2$. Then the following
conditions are equivalent.
\begin{enumerate}
\item[(a)] $T\in S_p$.
\item[(b)] $\{\|Te_n\|\}\in\ell^p$ for \emph{some} orthonormal basis $\{e_n\}$ in $H$.
\item[(c)] $\{\|Tf_n\|\}\in\ell^p$ for \emph{some} frame $\{f_n\}$ in $H$.
\end{enumerate}
Furthermore, we have
$$\|T\|^p_p=\inf\sum_{n=1}^\infty\|Te_n\|^p=\inf\sum_{n=1}^\infty\|f_n\|^{2-p}\|Tf_n\|^p
=\inf\sum_{n=1}^\infty\|Tf_n\|^p,$$
where the first infimum is taken over all orthonormal bases, the second infimum is taken
over all frames with lower frame bound greater than or equal to $1$, and the third infimum
is taken over all Parseval frames.
\label{17}
\end{thm}

\begin{proof}
Consider $S=T^*T$ and apply Theorem \ref{16} to the operator $S$. The desired result then
follows from the identity $\langle Sf_n,f_n\rangle=\|Tf_n\|^2$ and the fact that 
$T\in S_p$ if and only if $S\in S_{p/2}$.
\end{proof}

\begin{thm}
Suppose $0<p\le2$ and $T$ is a self-adjoint operator on $H$. Then the following
conditions are equivalent.
\begin{enumerate}
\item[(a)] $T\in S_p$.
\item[(b)] There exists \emph{some} orthonormal basis $\{e_n\}$ in $H$ such that
$$\sum_{k=1}^\infty\sum_{n=1}^\infty|\langle Te_n,e_k\rangle|^p<\infty.$$
\item[(c)] There exists \emph{some} frame $\{f_n\}$ in $H$ such that
$$\sum_{k=1}^\infty\sum_{n=1}^\infty|\langle Tf_n,f_k\rangle|^p<\infty.$$
\end{enumerate}
Moreover, we have
\begin{eqnarray*}
\|T\|_p^p&=&\inf\sum_{n=1}^\infty\sum_{k=1}^\infty|\langle Te_n,e_k\rangle|^p\\
&=&\inf\sum_{n=1}^\infty\sum_{k=1}^\infty\|f_n\|^{2-p}|\langle Tf_n,f_k\rangle|^p\\
&=&\inf\sum_{n=1}^\infty\sum_{k=1}^\infty|\langle Tf_n,f_k\rangle|^p,
\end{eqnarray*}
where the first infimum is taken over all orthonormal bases, the second infimum is taken
over all frames with lower frame bound at least $1$, and the third infimum is taken over 
all Parseval frames.
\label{18}
\end{thm}

\begin{proof}
If $T\in S_p$ is self-adjoint, then there exists an orthonormal set $\{\sigma_n\}$ 
such that
$$Tx=\sum_{n=1}^\infty\lambda_n\langle x,\sigma_n\rangle\,\sigma_n$$
for all $x\in H$, where $\{\lambda_n\}\in\ell^p$ is the nonzero eigenvalue sequence of 
$T$. Since each $\lambda_n$ is nonzero, we see that $Tx=0$ if and only if $x\perp\sigma_n$
for every $n$. Therefore, if $\{\sigma_n'\}$ is an orthonormal basis for $\ker(T)$, then
$\{e_n\}=:\{\sigma_n\}\cup\{\sigma_n'\}$ is an orthonormal basis for $H$. Moreover,
$$\sum_{n=1}^\infty\sum_{k=1}^\infty|\langle Te_n,e_k\rangle|^p=
\sum_{n=1}^\infty\sum_{k=1}^\infty|\langle T\sigma_n,\sigma_k\rangle|^p=
\sum_{n=1}^\infty|\lambda_n|^p<\infty.$$
This shows that condition (a) implies (b).

Since every orthonormal basis is a frame, it is trivial that condition (b) implies (c).
That (c) implies (a) follows from Theorem \ref{17} and part (b) of Lemma~\ref{11}.

If $\{f_n\}$ is a Parseval frame, then by Theorem~\ref{17} and the proof for part (b) 
of Lemma~\ref{11},
$$\|T\|_p^p\le\sum_{n=1}^\infty\|Tf_n\|^p\le\sum_{n=1}^\infty\sum_{k=1}^\infty
|\langle Tf_n,f_k\rangle|^p.$$
It follows that
$$\|T\|_p^p\le\inf\sum_{n=1}^\infty\sum_{k=1}^\infty|\langle Tf_n,f_k\rangle|^p.$$
Since every orthonormal basis is a Parseval frame, we clearly have
$$\inf\sum_{n=1}^\infty\sum_{k=1}^\infty|\langle Tf_n,f_k\rangle|^p\le
\inf\sum_{n=1}^\infty\sum_{k=1}^\infty|\langle Te_n,e_k\rangle|^p.$$
The inequality
$$\inf\sum_{n=1}^\infty\sum_{k=1}^\infty|\langle Te_n,e_k\rangle|^p\le\|T\|_p^p$$
follows from the first paragraph of this proof.

It follows from the proof of Lemma~\ref{11} that we always have
$$\sum_{n=1}^\infty\|f_n\|^{2-p}\sum_{k=1}^\infty|\langle Tf_n,f_k\rangle|^p
\ge\sum_{n=1}^\infty\|f_n\|^{2-p}\|Tf_n\|^p,$$
where $0<p\le2$ and $\{f_n\}$ is any frame with lower frame bound at least $1$.
Combining this with Theorem~\ref{17}, we see that
$$\|T\|_p^p=\inf\sum_{n=1}^\infty\sum_{k=1}^\infty\|f_n\|^{2-p}
|\langle Tf_n,f_k\rangle|^p,$$
where the infimum is taken over all frames with lower frame bound greater than or 
equal to $1$.
\end{proof}

We are not sure if the additional assumption that $T$ be self-adjoint is necessary in
Theorem~\ref{18} above. But we can show by an example that its proof will definitely not
work if no additional assumption is placed on $T$. To see this, fix any orthonormal 
basis $\{e_n\}$ and set
$$h_1=\sum_{n=1}^\infty\frac c{\sqrt n\,\log(n+1)}\,e_n,$$
where $c$ is a normalizing constant such that $\|h_1\|=1$. Expand $h_1$ to an orthonormal
basis $\{h_n\}$. Now define an operator $T$ on $H$ by
$$Tx=\sum_{n=1}^\infty\frac1{2^n}\langle x,h_n\rangle\,e_n,\qquad
T^*x=\sum_{n=1}^\infty\frac1{2^n}\langle x,e_n\rangle\,h_n.$$
It is easy to show that $T\in S_p$ for every $p>0$, but for $0<p<2$ we have
$$\sum_{n=1}^\infty\sum_{k=1}^\infty|\langle Te_n,e_k\rangle|^p=
\sum_{n=1}^\infty\frac1{2^{np}}\sum_{k=1}^\infty|\langle h_n,e_k\rangle|^p=\infty,$$
because in this case we have
$$\sum_{k=1}^\infty|\langle h_1,e_k\rangle|^p=c^p\sum_{k=1}^\infty
\frac1{[\sqrt k\,\log(k+1)]^p}=\infty.$$

It is now natural for us to ask whether the ``converse'' of the theorems above is true. 
We show that the answer is no when $p$ is not the upper end-point. The end-point case 
will be discussed in the next section. The next three propositions only require the
operator $T$ to be bounded, not necessarily in $S_p$.

\begin{prop}
Let $0<p<1$ and let $T$ be any nonzero operator on $H$. Then there exists a frame 
$\{f_n\}$ such that $\{\langle Tf_n,f_n\rangle\}\not\in\ell^p$.
\label{19}
\end{prop}

\begin{proof}
Since $T\not=0$, there exists a unit vector $h$ such that 
$\langle Th,h\rangle\not=0$. Fix such a vector $h$ and set
$$e_n'=\frac{h}{\sqrt n\log(n+1)},\qquad n\ge1.$$
For any $f\in H$ we have
$$\sum_{n=1}^\infty|\langle f,e_n'\rangle|^2=\sum_{n=1}^\infty\frac{|\langle f,h\rangle|^2}
{n[\log(n+1)]^2}\le\|f\|^2\sum_{n=1}^\infty\frac1{n[\log(n+1)]^2}.$$
Let $\{e_n\}$ be any orthonormal basis for $H$ and let $\{f_n\}=\{e_n\}\cup\{e_n'\}$.
Since the last series above converges and
$$\sum_{n=1}^\infty|\langle f,f_n\rangle|^2=\sum_{n=1}^\infty|\langle f,e_n\rangle|^2+
\sum_{n=1}^\infty|\langle f,e_n'\rangle|^2=
\|f\|^2+\sum_{n=1}^\infty|\langle f,e_n'\rangle|^2,$$
we conclude that $\{f_n\}$ is a frame for $H$.

On the other hand, the sequence $\{\langle Tf_n,f_n\rangle\}$ contains the subsequence
$\{\langle Te_n',e_n'\rangle\}$, which is not in $\ell^p$ for $0<p<1$. In fact,
$$\langle Te_n',e_n'\rangle=\frac{\langle Th,h\rangle}{n[\log(n+1)]^2}$$
for all $n\ge1$, which clearly shows that $\{\langle Te_n',e_n'\rangle\}\not\in\ell^p$
for $0<p<1$. This shows that $\{\langle Tf_n,f_n\rangle\}$ is not in $\ell^p$ and
completes the proof of the proposition.
\end{proof}

\begin{prop}
Let $0<p<2$ and let $T$ be any nonzero operator on $H$. There exists a frame 
$\{g_n\}$ such that $\{\|Tg_n\|\}\notin\ell^p$.
\label{20}
\end{prop}

\begin{proof}
Consider the operator $S=T^*T$ and use Proposition~\ref{19} to find a frame $\{f_n\}$
such that $\{\langle Sf_n,f_n\rangle\}$ is not in $\ell^{p/2}$. This is clearly the
same as $\{\|Tf_n\|\}\not\in\ell^p$.
\end{proof}

\begin{prop}
For any $0<p<2$ and any nonzero operator $T$ on $H$ there exists a frame $\{f_n\}$ such that
$$\sum_{k=1}^\infty\sum_{n=1}^\infty|\langle Tf_n,f_k\rangle|^p=\infty.$$
\label{21}
\end{prop}

\begin{proof}
This follows from Proposition~\ref{20} and part (b) of Lemma~\ref{11}.
\end{proof}

We can also obtain versions of these propositions in terms of orthonormal bases. Note that
the operators $T$ obtained in the next three propositions are in $S_p$, not just bounded.

\begin{prop}
Suppose $0<p<1$ and $\{e_n\}$ is any orthonormal basis for $H$. Then there exists a positive
operator $S\in S_p$ such that $\{\langle Se_n,e_n\rangle\}\not\in\ell^p$.
\label{22}
\end{prop}

\begin{proof}
Fix a nonzero operator $T\in S_p$ and use Proposition~\ref{19} to find a frame $\{f_n\}$ 
such that $\{\langle Tf_n,f_n\rangle\}\not\in\ell^p$. Let $A$ denote the operator from
Lemma~\ref{1} and consider the operator $S=A^*TA$. Then $S$ is a positive operator in
$S_p$ and $\{\langle Se_n,e_n\rangle\}=\{\langle Tf_n,f_n\rangle\}\not\in\ell^p$.
\end{proof}

\begin{prop}
Suppose $0<p<2$ and $\{e_n\}$ is any orthonormal basis for $H$. Then there exists a positive
operator $S\in S_p$ such that $\{\|Se_n\|\}\not\in\ell^p$.
\label{23}
\end{prop}

\begin{proof}
By Proposition~\ref{22}, there exists a positive operator $T\in S_{p/2}$ such that
$\{\langle Te_n,e_n\rangle\}\not\in\ell^{p/2}$. Let $S=\sqrt T$. Then $S$ is a positive
operator in $S_p$, $\|Se_n\|^2=\langle Te_n,e_n\rangle$, and $\{\|Se_n\|\}\not\in\ell^p$.
\end{proof}

\begin{prop}
For any $0<p<2$ and any orthonormal basis $\{e_n\}$ there exists a positive operator 
$S\in S_p$ such that 
$$\sum_{k=1}^\infty\sum_{n=1}^\infty|\langle Se_n,e_k\rangle|^p=\infty.$$
\label{24}
\end{prop}

\begin{proof}
This follows from Proposition~\ref{23} and part (b) of Lemma~\ref{11}. Note that the
operator constructed right before Proposition~\ref{19} is not positive.
\end{proof}

\section{The trace class and Hilbert-Schmidt class}

In this section we focus on two special classes of operators: the trace class $S_1$ and 
the Hilbert-Schmidt class $S_2$. It is well known that a bounded linear operator $T$ on $H$ 
is in $S_2$ if and only if $\sum\|Te_n\|^2<\infty$, where $\{e_n\}$ is any given 
orthonormal basis for $H$. Also, for $T\ge0$, $T\in S_1$ if and only if $\sum\langle 
Te_n,e_n\rangle<\infty$. We show that these results remain true when the orthonormal
basis $\{e_n\}$ is replaced by a frame.

\begin{thm}
Suppose $T$ is a positive operator on $H$. Then the following conditions are equivalent.
\begin{enumerate}
\item[(a)] $T$ is in the trace class $S_1$.
\item[(b)] $\{\langle Tf_n,f_n\rangle\}\in\ell^1$ for \emph{every} frame $\{f_n\}$.
\item[(c)] $\{\langle Tf_n,f_n\rangle\}\in\ell^1$ for \emph{some} frame $\{f_n\}$.
\end{enumerate}
\label{25}
\end{thm}

\begin{proof}
This follows from Theorems \ref{7} and \ref{16}.
\end{proof}

An equivalent version of Theorem \ref{25} above is the following.

\begin{thm}
Let $T$ be a bounded linear operator on $H$. Then the following conditions are equivalent.
\begin{enumerate}
\item[(a)] $T$ is in the Hilbert-Schmidt class $S_2$.
\item[(b)] $\{\|Tf_n\|\}\in\ell^2$ for \emph{every} frame $\{f_n\}$.
\item[(c)] $\{\|Tf_n\|\}\in\ell^2$ for \emph{some} frame $\{f_n\}$.
\end{enumerate}
\label{26}
\end{thm}

\begin{proof}
Note that $T$ is Hilbert-Schmidt if and only if $T^*T$ is trace class. Since
$\langle T^*Tf_n,f_n\rangle=\|Tf_n\|^2$, the desired result follows from Theorem~\ref{25}.
Alternatively, the desired result follows from Theorems~\ref{2} and \ref{17}.
\end{proof}

When $\{f_n\}$ is a frame, it is clear tha the condition $\{\|Tf_n\|\}\in\ell^2$ is
equivalent to the condition
$$\sum_{n=1}^\infty\sum_{k=1}^\infty|\langle Tf_n,f_k\rangle|^2<\infty.$$
Therefore, conditions (b) and (c) in Theorem~\ref{26} above can also be stated in terms of 
the double-indexed sequence $\{\langle Tf_n,f_k\rangle\}$.

\section{An application}

In this section, we consider a special class of frames in the Bergman space of the unit 
disk, namely, normalized reproducing kernels induced by sampling sequences. We use this
to obtain an integral condition for a bounded linear operator on the Bergman space to
belong to the Schatten class $S_p$.

Thus we let $A^2$ denote the space of analytic functions $f$ in the unit disk $\DD$ such that
$$\|f\|^2=\ind|f(z)|^2\,dA(z)<\infty,$$
where $dA$ is area measure on $\DD$ normalized so that $\DD$ has area $1$. As a closed subspace
of $L^2(\DD,dA)$, $A^2$ is a Hilbert space. In fact, $A^2$ is a reproducing Hilbert space whose
reproducing kernel is the well-known Bergman kernel
$$K_w(z)=K(z,w)=\frac1{(1-z\overline w)^2}.$$
For any $w\in\DD$ let $k_w$ denote the function in $A^2$ defined by
$$k_w(z)=\frac{K(z,w)}{\sqrt{K(w,w)}}=\frac{1-|w|^2}{(1-z\overline w)^2}.$$
Each $k_w$ is a unit vector in $A^2$, called the normalized reproducing kernel at $w$.

A sequence $\{w_n\}$ in $\DD$ is called a sampling sequence for the Bergman space $A^2$ if
there exists a positive constant $C$ such that
$$C^{-1}\|f\|^2\le\sum_{n=1}^\infty(1-|w_n|^2)^2|f(w_n)|^2\le C\|f\|^2$$
for all $f\in A^2$. This condition can be written as
$$C^{-1}\|f\|^2\le\sum_{n=1}^\infty|\langle f,k_{w_n}\rangle|^2\le C\|f\|^2.$$
Therefore, $\{w_n\}$ is a sampling sequence for the Bergman space if and only if the
sequence $\{k_{w_n}\}$ is a frame in $A^2$. See \cite{HKZ} for the theory
of Bergman spaces, including the notions of normalized reproducing kernels and
sampling sequences. Sampling sequences for the Bergman space are characterized in \cite{Seip}.

Some results obtained in the paper can be stated in terms of sampling sequences. As one 
particular example, we infer from Theorem~\ref{26} that if $\{w_n\}$ is a sampling 
sequence for the Bergman space, then a bounded linear operator $T$ on $A^2$ belongs to 
the Hilbert-Schmidt class $S_2$ if and only if it satisfies the condition
$$\sum_{n=1}^\infty\|Tk_{w_n}\|^2<\infty.$$
Equivalently, a positive operator $T$ on $A^2$ belongs to the trace class if and only if
$$\sum_{n=1}^\infty\langle Tk_{w_n},k_{w_n}\rangle<\infty.$$

\begin{lem}
Suppose $T$ is a bounded linear operator on $A^2$ and $0<p<\infty$. Then the function
$F(w)=\|TK_w\|^p$ is subharmonic in $\DD$.
\label{27}
\end{lem}

\begin{proof}
Without loss of generality we assume that $T\not=0$. It is then clear that the function 
$F$ has isolated zeros in $\DD$. Furthermore, away from the zeros of $F$, it follows from
$$F(w)=\langle TK_w,TK_w\rangle^{\frac p2}$$
that
$$\frac{\partial F}{\partial\overline w}=\frac p2\langle TK_w,TK_w\rangle^{\frac p2-1}
\langle TK_w^*,K_w\rangle,$$
where
$$K_w^*(z)=\frac{\partial}{\partial\overline w}\frac1{(1-z\overline w)^2}=
\frac{2z}{(1-z\overline w)^3}.$$
Differentiating one more time, we obtain
\begin{eqnarray*}
\frac{\partial^2F}{\partial w\partial\overline w}&=&\frac p2\left(\frac p2-1\right)
\langle TK_w,TK_w\rangle^{\frac p2-2}\langle TK_w^*,K_w\rangle\langle TK_w,TK_w^*\rangle\\
&&\ +\frac p2\langle TK_w,T_w\rangle^{\frac p2-1}\langle TK_w^*,TK_w^*\rangle\\
&=&\frac p2\!\left[\frac p2-1\right]\!\|TK_w\|^{p-4}|\langle TK_w^*,TK_w\rangle|^2
+\frac p2\|TK_w\|^{p-2}\|TK_w^*\|^2\\
&=&\frac p2\|TK_w\|^{p-4}\left[\left(\frac p2-1\right)|\langle TK_w^*,TK_w\rangle|^2
+\|TK_w\|^2\|TK_w^*\|^2\right]\\
&\ge&\frac p2\|TK_w\|^{p-4}\left[\|TK_w\|^2\|TK_w^*\|^2-|\langle TK_w^*,TK_w\rangle|^2\right]\\
&=&\frac p2\|TK_w\|^{p-4}\left[\|TK_w\|\|TK_w^*\|+|\langle TK_w^*,TK_w\rangle|\right]\cdot\\
&&\ \cdot\left[\|TK_w\|\|TK_w^*\|-|\langle TK_w^*,TK_w\rangle|\right]\\
&\ge&0.
\end{eqnarray*}
The last inequality above is a consequence of the Cauchy-Schwarz inequality.
\end{proof}

\begin{thm}
Suppose $T$ is a bounded linear operator on $A^2$ and
$$d\lambda(w)=\frac1{(1-|w|^2)^2}\,dA(w)$$
is the M\"obius invariant area measure on $\DD$. Then the condition
$$\ind\|Tk_w\|^p\,d\lambda(w)<\infty$$
is sufficient for $T\in S_p$ when $0<p\le2$ and it is necessary for $T\in S_p$ 
when $2\le p<\infty$. Consequently, $T$ is Hilbert-Schmidt if and only if
$$\ind\|Tk_w\|^2\,d\lambda(w)=\ind\|TK_w\|^2\,dA(w)<\infty.$$
\label{28}
\end{thm}

\begin{proof}
The result follows from Theorem 6.6 in \cite{Zhu}, because $T\in S_p$ if and only if
the positive operator $S=T^*T$ is in $S_{p/2}$,
$$\widetilde S(w)=\langle Sk_w,k_w\rangle=\|Tk_w\|^2$$
for all $w\in\DD$, and $0<p\le2$ if and only if $0<p/2\le1$.

The proof of Theorem 6.6 in \cite{Zhu} depends on the notion of the Berezin transform and
the spectral decomposition for positive operators. Here we give an independent proof in
the case $0<p\le2$ that is based on sampling sequences and subharmonicity.

Fix a sampling sequence $\{w_n\}\subset\DD$ for the Bergman space $A^2$ such that $\{w_n\}$ is
separated in the Bergman metric $\beta$, say, $\beta(w_i,w_j)>2r$ for some positive number $r$
and all $i\not=j$. See \cite{HKZ} for the existence of such a sequence. Let $D(w_n,r)=\{z\in\DD:
\beta(z,w_n)<r\}$ denote the Bergman metric ball centered at $w_n$ with radius $r$.

By Lemma \ref{27}, the function $w\mapsto\|TK_w\|^p$ is subharmonic. It follows from the 
proof of Proposition 4.13 in \cite{Zhu} that there exists a positive constant $C$, 
independent of $n$ and $T$, such that
$$\|TK_{w_n}\|^p\le\frac C{|D(w_n,r)|}\int_{D(w_n,r)}\|TK_w\|^p\,dA(w)$$
for all $n\ge1$, where $|D(w_n,r)|$ is the area of $D(w_n,r)$. By Proposition 4.5 in 
\cite{Zhu} and the remarks following it, we have
$$|D(w_n,r)|\sim(1-|w_n|^2)^2\sim(1-|w|^2)^2$$
for $w\in D(w_n,r)$. It follows that there exists another positive constant $C$, independent
of $n$ and $T$, such that
$$\|Tk_{w_n}\|^p\le C\int_{D(w_n,r)}\|Tk_w\|^p\,d\lambda(w)$$
for all $n\ge1$. Since the Bergman metric balls $D(w_n,r)$ are mutually disjoint, we have
$$\sum_{n=1}^\infty\|Tk_{w_n}\|^p\le C\ind\|Tk_w\|^p\,d\lambda(w).$$
The desired result now follows from Theorem \ref{17}.
\end{proof}

The ideas and results of this section clearly generalize to many other reproducing 
kernel Hilbert spaces, including weighted Bergman spaces on various domains and Fock spaces
on $\CC^n$.


\begin{thebibliography}{99}

\bibitem{ACP} Arias M.L., Corach G., and Pacheco M., Characterization of Bessel sequences, \textit{Extracta Mathematicae} 22 (2007), no. 1, 55--66.

\bibitem{Bal} Balazs P., Hilbert-Schmidt operators, frames classification, and best 
approximation by multipliers and algorithms, \textit{Int. J. Wavelets Multiresolut. Inf. Process}, 6 (2008), 315--330.

\bibitem{Chr} Christensen O., \textit{An Introduction to Frames and Riesz Basis}, Birkhauser, 2003.

\bibitem{GK} Goberg I.C. and Krein M.G., \textit{Introduction to the Theory of Linear
Nonselfadjoint Operators}, Translations of Mathematical Monographs Volume {\bf18},
American Mathematical Society, 1969.

\bibitem{HKZ} Hedenmalm H., Korenblum B., and Zhu K., \textit{Theory of Bergman Spaces},
Graduate Texts in Mathematics 199, Springer-Verlag, New York, 2000.

\bibitem{Koo} Koo Y.Y. and Lim J.K.  Schatten-class operators and frames, \textit{Questiones Mathematicase} 34 (2011), no. 2, 203--211.

\bibitem{Seip} Seip K., Beurling type density theorems in the unit disk, \textit{Invent. Math.}
113 (1993), 21-39.

\bibitem{Simon} Simon B., \textit{Trace Ideals and Their Applications}, Mathematical Surveys
and Monographs 120, Amer. Math. Soc., Providence, 2005.

\bibitem{Zhu} Zhu K., \textit{Operator Theory in Function Spaces}, 2nd edition, Mathematical Surveys and Monographs 138, Amer. Math. Soc., Providence, 2007.

\end{thebibliography}
\end{document}